\documentclass[11pt,a4paper,headinclude,footinclude,fleqn]{scrartcl}                 
\usepackage[T1]{fontenc}                   
\usepackage[utf8]{inputenc}                 
\usepackage[english]{babel}       
\usepackage{graphicx}                      %
\usepackage[font=small]{quoting}            %
\usepackage{caption}        
\usepackage{mathrsfs}          %
\usepackage[nochapters,beramono,eulermath,%
            pdfspacing,
            listings,
            ]{classicthesis}                %
\usepackage{arsclassica}                    
\usepackage[top=.8in,bottom=1in,left=1in,right=1in]{geometry}
\usepackage[english]{babel}
\usepackage{verbatim}
\usepackage{color}
\usepackage{picture}
\usepackage{amsmath,amsthm,amsfonts,amssymb}
\usepackage{comment}
\usepackage{mathtools}

\usepackage{titling} 
\usepackage[hyperpageref]{backref}

\newtheorem{thm}{Theorem}[section]
\newtheorem{lem}[thm]{Lemma}

\newtheorem{prop}[thm]{Proposition}
\newtheorem{example}{Example}

\theoremstyle{definition}
\newtheorem{rem}[thm]{Remark}
\theoremstyle{remark}

\newcommand{\ds}{\displaystyle}

{\left\{\begin{array}{@{}l@{}}}{\end{array}\right.}
\patchcmd{\abstract}{\scshape\abstractname}{\textbf{\abstractname}}{}{}
\makeatletter 
\def\@makefnmark{} 
\makeatother 
\title{A sharp weighted anisotropic Poincar\'e inequality \\ for convex domains}
\author{Francesco Della Pietra, %
Nunzia Gavitone, %
Gianpaolo Piscitelli \\[.3cm]%
{\em\scriptsize Universit\`a degli studi di Napoli Federico II, Dipartimento di Matematica e Applicazioni ``R. Caccioppoli''} \\ 
{\em\scriptsize Via Cintia, Monte S. Angelo - 80126 Napoli, Italia.}\thanks{Email: f.dellapietra@unina.it, nunzia.gavitone@unina.it, gianpaolo.piscitelli@unina.it
}
}

\begin{document}

\maketitle


\begin{abstract}
We prove an optimal lower bound for the best constant in a class of weighted anisotropic Poincar\'e inequalities.\end{abstract}

\section{Introduction}
\label{intro}
In this paper we prove a sharp lower bound for the optimal constant 
$\mu_{p,\mathcal H,\omega} (\Omega)$ in the Poincar\'e-type inequality
\begin{equation*}
\inf_{t\in \mathbb R} 
\|u-t\|_{L^p_{\omega}(\Omega)}\le \frac{1}{[\mu_{p,\mathcal H,\omega}(\Omega)]^{\frac1p}}\|\mathcal H(\nabla u)\|_{L^p_{\omega}(\Omega)},
\end{equation*}
with $1<p<+\infty$, $\Omega$ is a bounded convex domain of $\mathbb R^{n}$, $\mathcal H\in \mathscr H(\mathbb R^{n})$, where $\mathscr H(\mathbb R^{n})$ is the set of lower semicontinuous functions, positive in $\mathbb R^{n}\setminus\{0\}$ and positively $1$-homogeneous; moreover, let $\omega$ be a $\log$-concave function.

If $\mathcal H$ is the Euclidean norm of $\mathbb R^{n}$ and $\omega=1$, then $\mu_{p}(\Omega)=\mu_{p,\mathcal E,\omega}(\Omega)$ is the first nontrivial eigenvalue of the Neumann $p-$Laplacian:
\begin{equation*}
\left\{
\begin{array}{ll}
-\Delta_{p} u=\mu_{p}|u|^{p-2}u &\text{in }\Omega,\\[.2cm]
|\nabla u|^{p-2} \frac{\partial u}{\partial \nu}=0 &\text{on }\partial\Omega.
\end{array}
\right.
\end{equation*}
Then, for a convex set $\Omega$ it holds that
\begin{equation*}
\mu_{p}(\Omega) \ge \left(\frac{\pi_p}{D_\mathcal E(\Omega)}\right)^p,
\end{equation*}
where
\[
\pi_p=2\int_0^{+\infty}\frac{1}{1+\frac{1}{p-1}s^p}ds=2\pi\frac{(p-1)^{\frac1p}}{p\sin\frac{\pi}{p}},\qquad D_{\mathcal E}(\Omega)\text{ Euclidean diameter of }\Omega.
\]
This estimate, proved in the case $p=2$ in \cite{pw} (see also \cite{bebe}), has been generalized the case $p\ne 2$ in \cite{ad,ent,fnt,v} and for $p\to \infty$ in \cite{eknt,rs}. Moreover the constant $\left(\frac{\pi_p}{D_\mathcal E(\Omega)}\right)^p$ is the optimal constant of the one-dimensional Poincar\'e-Wirtinger inequality, with $\omega=1$, on a segment of length $D_{\mathcal E}(\Omega)$. When $p=2$ and $\omega=1$, in \cite{bcdl} an extension of the estimate in the class of suitable non-convex domains has been proved.

The aim of the paper is to prove an analogous sharp lower bound for $\mu_{p,\mathcal H,\omega}(\Omega)$, in a general anisotropic case. More precisely, our main result is:
\begin{thm}
\label{main}
Let $\mathcal H\in\mathscr H(\mathbb R^{n})$, $\mathcal H^{o}$ be its polar function. Let us consider a bounded convex domain $\Omega\subset\mathbb R^n$, $1<p<\infty$, and take a positive $\log$-concave function $\omega$ defined in $\Omega$. Then, given
\[
\mu_{p,\mathcal H,\omega}(\Omega)=\inf_{\substack{u\in W^{1,\infty}(\Omega) \\ \int_\Omega |u|^{p-2}u\omega\,dx=0}} \dfrac{\displaystyle\int_\Omega \mathcal H(\nabla u)^p\omega\,dx}{\displaystyle\int_\Omega |u|^p\omega\,dx},
\]
it holds that
\begin{equation}
\label{pwineq}
\mu_{p,\mathcal H,\omega}(\Omega) \ge \left(\frac{\pi_p}{D_\mathcal H(\Omega)}\right)^p,
\end{equation}
where $D_\mathcal H(\Omega)=\sup_{x,y\in\Omega}\mathcal H^{o}(y-x)$.
\end{thm}

 This result has been proved in the case $p=2$ and $\omega=1$, when $\mathcal H$ is a strongly convex, smooth norm of $\mathbb R^{n}$ in \cite{wx} with a completely different method than the one presented here. 

In Section 2 below we give the precise definition of $\mathcal H^{o}$ and give some details on the set $\mathscr H(\mathbb R^{n})$. In Section 3 we give the proof of the main result.



\section{Notation and preliminaries}
\label{notation}
A function
\[
\xi\in \mathbb R^{n}\mapsto \mathcal H(\xi)\in [0,+\infty[
\] 
belongs to the set $\mathscr H(\mathbb R^{n})$
if it verifies the following assumptions:
\begin{enumerate}
\item \label{p1} $\mathcal H$ is positively 1-homogeneous, that is 
\begin{equation*}
\text{if }\xi\in \mathbb R^{n}\text{ and }t\ge 0,\text{ then } \mathcal H(t\xi)=t\mathcal H(\xi);
\end{equation*} 
\item \label{p2} if $\xi\in\mathbb R^{n}\setminus\{0\}$, then $\mathcal H(\xi)>0$;
\item \label{p3} $\mathcal H$ is lower semi-continuous.
\end{enumerate}
If $\mathcal H\in \mathscr H(\mathbb R^{n})$, properties \eqref{p1}, \eqref{p2}, \eqref{p3} give that there exists a positive constant $a$ such that 
\begin{equation*}
\label{eq:lin}
a|\xi| \le \mathcal H(\xi),\quad \xi \in \mathbb R^{n}.
\end{equation*} 

The polar function $\mathcal H^o\colon\mathbb R^n \rightarrow [0,+\infty[$ 
of $\mathcal H \in \mathscr H(\mathbb R^{n})$ is defined as
\begin{equation*}
\mathcal H^o(\eta)=\sup_{\xi \ne 0} \frac{\langle \xi, \eta\rangle}{\mathcal H(\xi)}. 
\end{equation*}
The function $\mathcal H^o$ belongs to $\mathscr H(\mathbb R^{n})$. Moreover it is convex on $\mathbb R^{n}$, and then continuous.
If $\mathcal H$ is convex, it holds that 
\[
\mathcal H(\xi)=(\mathcal H^{o})^{o}(\xi)=\sup_{\eta \ne 0} \frac{\langle \xi, \eta\rangle}{\mathcal H^{o}(\eta)}.
\] 
If $\mathcal H$ is convex and $\mathcal H(\xi)=\mathcal H(-\xi)$ for all $\xi \in \mathbb R^{n}$, then $\mathcal H$ is a norm on $\mathbb R^{n}$, and the same holds for $\mathcal H^o$.

We recall that if $\mathcal H$ is a smooth norm of $\mathbb R^{n}$ such that $\nabla^{2}(\mathcal H^{2})$ is positive definite on $\mathbb R^{n}\setminus\{0\}$, then $\mathcal H$ is called a Finsler norm on $\mathbb R^{n}$. 

If $\mathcal H\in \mathscr H(\mathbb R^{n})$, by definition we have
\begin{equation}
\label{riduzione}
 \langle \xi, \eta\rangle  \le \mathcal H(\xi) \mathcal H^{o}(\eta), \qquad \forall \xi, \eta \in \mathbb R^{n}.
\end{equation}

\begin{rem}
Let $\mathcal H\in \mathscr H(\mathbb R^{n})$, and consider the convex envelope of $\mathcal H$, that is the largest convex function $\overline{\mathcal H}$ such that $\overline{\mathcal H}\le \mathcal H$. It holds that $\overline{\mathcal H}$ and $\mathcal H$ have the same polar function: 
\[
(\overline{\mathcal H})^{o}= \mathcal H^{o}\quad \text{in }\mathbb R^{n}.
\]
Indeed, being $\overline{\mathcal H}\le \mathcal H$, by definition it holds that $(\overline{\mathcal H})^{o} \ge \mathcal H^{o}$. To show the reverse inequality, it is enough to prove that $(\mathcal H^{o})^{o} \le \mathcal H$. Then, being $\overline{\mathcal H}$ the convex envelope of $\mathcal H$, it must be $(\mathcal H^{o})^{o} \le \overline{\mathcal H}$, that implies $(\overline{\mathcal H})^{o} \le \mathcal H^{o}$. Denoting by $G(x)=(\mathcal H^{o})^{o}(x)$, for any $x$ there exists $\overline v_{x}$ such that
\[
G(x)=\frac{\langle x,\overline v_{x}\rangle}{\mathcal H^{o}(\overline v_{x})},\quad \text{and}\quad \langle x,\overline v_{x}\rangle \le \mathcal H^{o}(\overline v_{x})\mathcal H(x), \quad\text{that implies}\quad G(x)\le \mathcal H(x).
\]
\end{rem}

Let $\mathcal H\in \mathscr H(\mathbb R^{n})$, and consider a bounded convex domain $\Omega$ of $\mathbb R^{n}$. Throughout the paper $D_\mathcal H(\Omega)\in ]0,+\infty[$ will be
\begin{equation*}
D_\mathcal H(\Omega)=\sup_{x,y\in\Omega} \mathcal H^{o}(y-x).
\end{equation*}
We explicitly observe that since $\mathcal H^{o}$ is not necessarily even, in general $\mathcal H^{o}(y-x)\ne \mathcal H^{o}(x-y)$. When $\mathcal H$ is a norm, then $D_\mathcal H(\Omega)$ is the so called anisotropic diameter of $\Omega$ with respect to $\mathcal H^{o}$. In particular, if $\mathcal H=\mathcal E$ is the Euclidean norm in $\mathbb R^{n}$, then $\mathcal E^{o}=\mathcal E$ and  $D_{\mathcal E}(\Omega)$ is the standard Euclidean diameter of $\Omega$.  We refer the reader, for example, to \cite{cs,ffk} for remarkable examples of convex not even functions in $\mathscr H(\mathbb R^{n})$. On the other hand, in \cite{vs} some results on isoperimetric and optimal Hardy-Sobolev inequalities for a general function $\mathcal H\in \mathscr H(\mathbb R^{n})$ have been proved, by using a generalizazion of the so called convex symmetrization introduced in \cite{aflt} (see also \cite{dgmana,dgpota,dgmaan}).
\begin{rem}
In general $\mathcal H$ and $\mathcal H^{o}$ are not rotational invariant. Anyway, if $A\in SO(n)$, defining 
\begin{equation}
\label{HA}
\mathcal H_{A}(x) =\mathcal H(Ax),
\end{equation}
and being $A^{T}=A^{-1}$, then $\mathcal H_{A}\in \mathscr H(\mathbb R^{n})$ and
\[
(\mathcal H_{A})^{o}(\xi)=\sup_{x\in \mathbb R^{n}\setminus\{0\}}\frac{\langle x,\xi\rangle}{\mathcal H_{A}(x)}=
\sup_{y\in \mathbb R^{n}\setminus\{0\}}\frac{\langle A^{T}y,\xi\rangle}{\mathcal H(y)}=
\sup_{y\in \mathbb R^{n}\setminus\{0\}}\frac{\langle y, A\xi\rangle}{\mathcal H(y)}= (\mathcal H^{o})_{A}(\xi).
\]
Moreover,
\begin{equation}
\label{rotazD}
D_{\mathcal H_A}(A^{T}\Omega)=\sup_{x,y\in A^{T}\Omega} (\mathcal H^{o})_{A}(y-x)=
\sup_{\bar x,\bar y\in \Omega} \mathcal H^{o}(\bar y-\bar x) =D_{\mathcal H}(\Omega).
\end{equation}
\end{rem}







\section{Proof of the Payne-Weinberger inequality}
In this section we state and prove Theorem \ref{main}. To this aim, the following Wirtinger-type inequality, contained in \cite{fnt} is needed.
\begin{prop}
\label{propuny}
Let $f$ be a positive $log$-concave function defined on $[0,L]$ and $p>1$, then 
\begin{equation*}
\inf
 \left\{
\frac{\displaystyle\int_0^{L} |u'|^p f \,dx}{\displaystyle\int_0^{L} |u|^p f \,dx},\; u \in W^{1,p}(0,L),\; \int_0^L |u|^{p-2}u f dx=0 \right\}
 \ge \frac{\pi_p^p}{L^p}.
\end{equation*}
\end{prop}

The proof of the main result is based on a slicing method introduced in \cite{pw} in the Laplacian case. The key ingredient is the following Lemma. For a proof, we refer the reader, for example, to \cite{pw,bebe,fnt}.
\begin{lem}
\label{keylem}
Let $\Omega$ be a convex set in $\mathbb R^n$ having (Euclidean) diameter $D_{\mathcal E}(\Omega)$, let $\omega$ be a positive log-concave function on $\Omega$, and let $u$ be any function such that $\int_{\Omega} |u|^{p-2}u \omega \, dx=0$. Then, for all positive $\varepsilon$, there exists a decomposition of the set $\Omega$ in mutually disjoint convex sets $\Omega_i$ $(i=1,\ldots,k)$ such that
\begin{gather*}
\bigcup_{i=1}^{k} \overline{\Omega}= \overline{\Omega} \\
\int_{\Omega_i} |u|^{p-2}u \omega \, dx=0
\end{gather*}
and for each $i$ there exists a rectangular system of coordinates such that
\[
\Omega_i \subset \{(x_1,\ldots, x_n) \in \mathbb R^n \colon 0\le x_1 \le d_i, |x_l| \le \varepsilon,\, l=2,\ldots,n\},
\]
where $d_i \le D_{\mathcal E}(\Omega)$, $i=1, \ldots,k$. 
\end{lem}

\textbf{Proof of Theorem \ref{main}}.
By density, it is sufficient to consider a smooth function $u$ with uniformly continuous first derivatives and $\int_{\Omega}|u|^{p-2}u \omega \, dx=0$.

Hence, we can decompose the set $\Omega$ in $k$ convex domains $\Omega_i$ as in Lemma \ref{keylem}. 
In order to prove \eqref{pwineq}, we will show that for any $i\in\{1,\ldots,k\}$ it holds that
\begin{equation}
\label{passo}
\int_{\Omega_{i}} H^{p}(\nabla u)\omega\,dx \ge \frac{\pi_{p}^{p}}{D_{\mathcal H}(\Omega)^{p}}\int_{\Omega_{i}}|u|^{p}\omega\,dx.
\end{equation}
By Lemma \ref{keylem}, for each fixed $i\in\{1,\ldots,k\}$, there exists a rotation $A_{i}\in SO(n)$ such that
\begin{equation*}
A_{i}\Omega_i \subset \{(x_1,\ldots, x_n) \in \mathbb R^n \colon 0\le x_1 \le d_i,\; |x_l| \le \varepsilon,\; l=2,\ldots,n\}.
\end{equation*}
By changing the variable $y=A_{i}x$, recalling the notation \eqref{HA} and using   \eqref{rotazD} it holds that
\[
\int_{\Omega_{i}} \mathcal H^{p}(\nabla u(x))\,\omega(x)\,dx=
\int_{A_{i}\Omega_{i}} \mathcal H_{A_{i}^{T}}(\nabla u(A_{i}^{T}y))^{p}\,\omega(A_{i}^{T}y)\,dy;\qquad D_{\mathcal H}(\Omega)=D_{\mathcal H_{A_{i}^{T}}}(A_{i}\Omega).
\]
We deduce that it is not restrictive to suppose that for any $i\in\{1,\ldots,n\}$ $A_{i}$ is the identity matrix, and the decomposition holds with respect to the $x_{1}-$axis.

Now we may argue as in \cite{fnt}. For any $t \in [0,d_i]$ let us denote by $v(t)=u(t,0,\ldots,0),$ and $f_i(t)=g_i(t) \omega(t,0,\ldots,0)$, where $g_i(t)$ will be the $(n-1)$ volume of the intersection of $\Omega_i$ with the hyperplane $x_1=t$. By Brunn-Minkowski inequality $g_i$, and then $f_i$, is a log-concave function in $[0,d_i]$.
Since $u, u_{x_1}$ and $\omega$ are uniformly continuous in $\Omega$ there exists a  modulus of continuity $\eta(\cdot)$ with $\eta(\varepsilon) \searrow 0$ for $\varepsilon \to 0$, indipendent of the decomposition of $\Omega$ and such that
\begin{equation*}
 \left| \int_{\Omega_i} |u_{x_1}|^p\omega\, dx - \int_0^{d_i} |v'|^p f_i \,dt\right| \le \eta(\varepsilon) |\Omega_i|, \quad\qquad 
\left| \int_{\Omega_i} |u|^p\omega\, dx - \int_0^{d_i} |v|^p f_i \,dt\right| \le \eta(\varepsilon) |\Omega_i|, 
\end{equation*}
and
\begin{equation*}
\left|\int_0^{d_i} |v|^{p-2}v f_i \,dt\right|\le \eta(\varepsilon) |\Omega_i|.
\end{equation*}
Now, by property \eqref{riduzione} we deduce that for any vector $\eta\in\mathbb R^{n}$
\[
|\langle \nabla u,\eta\rangle| \le \mathcal H(\nabla u) \max\{ \mathcal H^o(\eta),\mathcal H^o(-\eta)\}.
\]
Then choosing $\eta=e_{1}$ and denoting by $M=\max\{ \mathcal H^o(e_{1}),\mathcal H^o(-e_{1})\}$, Proposition \ref{propuny} gives
\begin{multline*}
\int_{\Omega_i}\mathcal H^p(\nabla u) \omega\, dx \ge \frac{1}{M^p} \int_{\Omega_i} |u_{x_1}|^p\omega\, dx \ge  \frac{1}{M^p}  \int_0^{d_i} |v'|^p f_i \,dt-\frac{\eta(\varepsilon) |\Omega_i|}{M^p}\\ \ge   \frac{\pi_p}{d_i^p M^p}  \int_0^{d_i} |v|^p f_i \,dt +C\eta(\varepsilon) |\Omega_i| \ge  \frac{\pi_p^p}{d_i^pM^p}  \int_{\Omega_i} |u|^p \omega \,dx +C\eta(\varepsilon) |\Omega_i|, 
\end{multline*}
where $C$ is a constant which does not depend on $\varepsilon$.
Being $d_i\le D_\mathcal E(\Omega)$, and then $d_i M\le D_\mathcal H(\Omega)$, by letting $\varepsilon$ to zero we get \eqref{passo}. Hence, by summing over $i$ we get the thesis.

\begin{rem}
In order to prove an estimate for $\mu_{p,\mathcal H,\omega}$, we could use directly property \eqref{riduzione} with $v= \frac{\nabla u}{|\nabla u|}$, and the Payne-Weinberger inequality in the Euclidean case, obtaining that
\[
\int_{\Omega}\mathcal H^p(\nabla u) \omega\, dx \ge \int_{\Omega} \frac{|\nabla u|^p}{\mathcal H^o(v)^p}\omega\, dx \ge  \frac{\pi_p^p}{D_{\mathcal E}(\Omega)^p \mathcal H^o(v_m)^p}
\int_{\Omega} |u|^p \omega \,dx, 
\]
where $\ds \mathcal H^o(v_m)=\max_{|\nu|=1}\mathcal H^o(\nu)$. However, we have a worst estimate than \eqref{pwineq} because  $D_{\mathcal E}(\Omega) \cdot \mathcal H^o(v_m)$ is, in general, strictly larger than $D_\mathcal H(\Omega)$, as shown in the following example.
\end{rem}

\vspace{.3cm}

\begin{example}
Let $\mathcal H(x,y)=\sqrt{a^2 x^2+b^2y^2}$, with $a<b$. Then $\mathcal H$ is a  even, smooth norm with  $\mathcal H^o(x,y)=\sqrt{ \dfrac{x^2}{a^2}+\dfrac{y^2}{b^2}}$ and the Wulff shapes $\{\mathcal H^{o}(x,y)< R\}$, $R>0$, are ellipses. Clearly we have:
\[
D_{\mathcal E}(\Omega)=2b \quad \text{ and }\quad D_\mathcal H(\Omega)=2
\]
Let us compute  $\mathcal H^o(v_m)$. We have: 
\[
\max_{|v|=1}\mathcal H^o(v)= \max_{\vartheta \in [0,2 \pi]}\sqrt{\frac{(\cos \vartheta)^2}{a^2}+\frac{(\sin \vartheta)^2}{b^2}} = \mathcal H^o(0,\pm 1)= \frac{1}{a}.
\]
Then $D_{\mathcal E}(\Omega)\cdot \mathcal H^o(v_m)=2 \dfrac{b}{a} > 2$. 
\end{example}

\end{document}